\documentclass[12pt,psamsfonts]{amsart}

\usepackage{amsmath,amssymb}
\usepackage{graphicx}
\usepackage{color}
\usepackage{psfrag} 
\usepackage{overpic} 

\newtheorem{theorem}{Theorem}

\newtheorem {remark}[theorem]{Remark}
\newtheorem{definition}[theorem]{Definition}

\title[Cyclicity for Quadratic three-dimensional systems]{Lower Bounds for the Cyclicity of centers of Quadratic three-dimensional systems}
\author[L.F. Gouveia, L. Queiroz]{}


\begin{document}
 \maketitle

\centerline{\scshape  Luiz F. S. Gouveia,  \; Lucas Queiroz}
\medskip

{\footnotesize \centerline{Universidade Estadual Paulista (UNESP), Instituto de Bioci\^encias Letras e Ci\^encias Exatas,} \centerline{R. Cristov\~ao Colombo, 2265, 15.054-000, S. J. Rio Preto, SP, Brasil }
\centerline{\email{fernando.gouveia@unesp.br} and \email{lucas.queiroz@unesp.br}}}

\medskip

\bigskip

\begin{quote}{\normalfont\fontsize{8}{10}\selectfont
{\bfseries Abstract.}
We consider quadratic three-dimensional differential systems having a Hopf singular point. We study the cyclicity when the singular point is a center on the center manifold using higher order developments of the Lyapunov constants. As a result, we make a chart of the cyclicity by establishing the lower bounds for several known systems in the literature, among them the R{\"o}ssler, Lorenz and Moon-Rand systems. Moreover, we obtain an example of a jerk system for which is possible to bifurcate 12 limit-cycles from the center, which is a new lower bound for three-dimensional quadratic systems.
\par}
\end{quote}

\section{Introduction}
The problem of determining and locating the limit-cycles for a given planar polynomial differential system is known as the \emph{Hilbert's sixteenth problem} and dates back to the year 1900. In the last century, several researchers have approached this problem and made great advancements, which are described in great detail in the surveys \cite{IlyashenkoSurvey, JibinLi} by Ilyashenko and Jibin Li respectively. However, Hilbert's sixteenth problem has not been completely solved. There are several simpler versions of the problem, among them we can highlight the problem of determining the quantity $M(n)$ of small amplitude limit-cycles which can bifurcate from an elementary center or focus of a polynomial vector field of degree $n$ (see, for instance, \cite{Zoladek1}). A singular point of a planar polynomial vector field is an \emph{elementary center or focus} if the eigenvalues of its Jacobian matrix evaluated at the singular point are purely imaginary. Any differential system associated to a vector field having an elementary center or focus can be written, after the proper change of variables and time rescaling, in the form:
\begin{equation*}\label{eqPlanar1}
\left\{\begin{array}{lcr}
\dot{x}=-y+P(x,y),\\
\dot{y}=x+Q(x,y),
\end{array}\right.
\end{equation*}
where $P$ and $Q$ are polynomials with no linear or constant terms.

Several techniques were developed to study the cyclicity, that is, the potential to bifurcate limit-cycles, of an elementary center or focus. Among those, one of the most powerful is the computation of \emph{Lyapunov constants} which are also related to another important problem in the qualitative theory of differential equations, the so-called \emph{Center Problem} \cite{ArtDumLli2006, Romanovski}. Christopher in \cite{Christopher} proved that by performing an analysis on the linear part of the Lyapunov constants in its power series expansion with respect to the perturbation parameters, it is possible to estimate the cyclicity of a center. Using this idea, Torregrosa and Liang \cite{TorrePara} proposed the \emph{Parallelization method} which allow computations for the linear parts to be made separately, reducing computational time. This method was further explored in \cite{TorreGouveia} where the authors improved the previous lower bounds of the cyclicity for centers of some polynomial vector fields, using higher order developments of the Lyapunov constants with respect to the perturbation parameters.

A natural question that arises is if it is possible to extend the cyclicity problem and the methods to approach it to three-dimensional differential systems. And in fact, the answer is positive. We consider analytical vector fields in $\mathbb{R}^3$ having a \emph{Hopf singular point}, that is, a singular point for which the Jacobian matrix has a pair of purely imaginary eigenvalues and a non-zero real eigenvalue. The differential system associated to such vector fields can be put in the following canonical form, by means of a linear change of variables and time rescaling:
\begin{equation}\label{eq1}
\left\{\begin{array}{lcr}
\dot{x}=-y+P(x,y,z),\\
\dot{y}=x+Q(x,y,z),\\
\dot{z}=-\lambda z+R(x,y,z),
\end{array}\right.
\end{equation}
where $P,Q,R$ are polynomials with no linear nor constant terms and $\lambda\neq 0$. For vector fields \eqref{eq1} there exists an invariant bidimensional $C^r$-manifold $W^c$ tangent to the $xy$-plane at the origin for every $r\geqslant 3$. This result is known as the \emph{Center Manifold Theorem}, and its proof, along with a more detailed study can be found in \cite{Kelley, Sijbrand}. 

The restriction of \eqref{eq1} to a center manifold $W^c$ is a bidimensional differential system which has an elementary center or focus and we can investigate the bifurcation of limit-cycles. Although the above theorem insures that the invariant manifold exists, neither its analyticity nor its uniqueness is guaranteed. However, it is known that the flow of the restriction to any two $C^r$ center manifolds are $C^{r-1}$-conjugated \cite{Aulbach, Burchard}. 

Since the center manifold is not necessarily unique nor analytic, computing a parametrization for $W^c$ and then applying the planar techniques is not optimal. Fortunately, the computation of the Lyapunov constants can be made even without knowing a parametrization for any center manifold \cite{Adam}. Furthermore, in \cite{Garcia} the authors proved that it is possible to study cyclicity for the Hopf singularity on center manifolds via the Lyapunov constants in the same way proposed by Christopher. The authors of \cite{Ivan} used this approach to give new lower bounds for quadratic, cubic, quartic and quintic three-dimensional systems. The respective lower bounds are 11, 31, 54 and 92 limit-cycles.

Our objective is to make a chart of the cyclicity for the known quadratic systems \eqref{eq1} in the literature which have a center on the center manifold. In order to attain this goal, we use the high order developments of the Lyapunov constants and full quadratic perturbations. Moreover, we exhibit an example of a quadratic three-dimensional system for which is possible to bifurcate 12 limit-cycles from the center on the center manifold (Theorem \ref{TeoJerk12}). This is so far the best known lower bound in the literature. 

The main results of this work are summed up in the following theorems:

\begin{theorem}\label{TeoLorenzRosslerMoonRand}
There exist quadratic perturbations of the R{\"o}ssler, Lorenz and Moon-Rand systems for which there are 4 limit-cycles bifurcating from the origin.
\end{theorem}

\begin{theorem}\label{TeoJerk12}
The origin of the following quadratic jerk system
\begin{equation}\label{eqJerk12}
\left\{\begin{array}{lll}
		\dot{x}&=&-y+2x^2+4xz-4y^2+8yz-10z^2,\vspace{0.25cm} \\ 
		\dot{y}&=&x+2x^2+4xz-4y^2+8yz-10z^2,\vspace{0.25cm} \\ 
		\dot{z}&=&-z+2x^2+4xz-4y^2+8yz-10z^2,
	\end{array}
	\right.
\end{equation}
is a center on the center manifold and unfolds 12 limit-cycles under quadratic perturbations.
\end{theorem}

\begin{theorem}\label{TeoGineValls}
There exist parameter values for which the origin of the system
\begin{equation}\label{Valls}
	\left\{\begin{array}{lcr}
		\dot{x}=y,\\
		\dot{y}=-x+a_1x^2+a_2xy+a_3xz+a_4y^2+a_5yz+a_6z^2,\\
		\dot{z}=-z+c_1x^2+c_2xy+c_3y^2,	
	\end{array}\right.	
\end{equation}
is a center on the center manifold and unfolds 10 limit-cycles under quadratic perturbations.
\end{theorem}

\begin{theorem}\label{TeoAdam}
There exist parameter values for which the origin of the system
\begin{equation}\left\{\begin{array}{lll}\label{Adam}
		\dot{x}&=&-y+ax^2+ay^2+cxz+dyz, \vspace{0.25cm} \\ 
		\dot{y}&=&x+bx^2+by^2+exz+fyz, \vspace{0.25cm} \\
		\dot{z}&=&-z+Sx^2+Sy^2+Txz+Uyz.
	\end{array}
	\right.
\end{equation}
is a center on the center manifold and unfolds 9 limit-cycles under quadratic perturbations.
\end{theorem}

The structure of this paper is as follows: In Section 2, we exhibit some fundamental concepts and results necessary for the development of our investigation. In Section 3, we study celebrated systems in the literature proving Theorem \ref{TeoLorenzRosslerMoonRand} which proposes new lower bounds for the cyclicity of the R{\"o}ssler \cite{MalasomaRossler}, Lorenz \cite{Lorenz}, and Moon-Rand \cite{MoonRand} systems. Section 4 is dedicated to study the jerk systems, where we prove Theorem \ref{TeoJerk12} obtaining a new lower bound for the cyclicity of three-dimensional quadratic systems of 12 limit-cycles. In sections 5, we study the cyclicity of the quadratic systems \eqref{Valls} and \eqref{Adam} considered in the papers \cite{Gine} and \cite{Adam} respectively, proving Theorems \ref{TeoGineValls} and \ref{TeoAdam}.

\section{Lyapunov constants and Cyclicity}

One of the most useful tools to study monodromic singular points is \emph{Poincar\'e map} or the \emph{First Return map}. For planar systems, it is widely known how to define those maps and its properties. In \cite{Garcia} the authors show how to extend those concepts for system \eqref{eq1}. We state the main results here, and encourage the reader to refer to \cite{Queiroz, GarciaJacobi, Garcia} for more details and proofs. 

Introducing the change of variables $x=\rho\cos\theta$, $y=\rho\sin\theta$ and $z=\rho\omega$, we can describe the solution curves of \eqref{eq1} by the following equations:
\begin{equation}\label{eqRW}
\left\{\begin{array}{lcr}
\dfrac{d\rho}{d\theta}=\dfrac{P\cos\theta+Q\sin\theta}{1+\rho[\tilde{Q}\cos\theta-\tilde{P}\sin\theta]},\\
\dfrac{d\omega}{d\theta}=\dfrac{-\mu\omega+\rho[\tilde{R}-\omega(\tilde{P}\cos\theta+\tilde{Q}\sin\theta)]}{1+\rho[\tilde{Q}\cos\theta-\tilde{P}\sin\theta]}.
\end{array}\right.
\end{equation}
For each $(\rho_0,\omega_0)$ with sufficiently small $\|(\rho_0,\omega_0)\|$, let $\varphi(\theta,\rho_0,\omega_0)=(\rho(\theta,\rho_0,\omega_0),\omega(\theta,\rho_0,\omega_0))$ be the solution of system \eqref{eqRW} with initial conditions $(\rho(0,\rho_0,\omega_0),\omega(0,\rho_0,\omega_0))=(\rho_0,\omega_0)$. We then define:
\begin{definition}For system \eqref{eq1}, the map $$d(\rho_0,\omega_0)=\varphi(2\pi,\rho_0,\omega_0)-(\rho_0,\omega_0)$$ is called \emph{displacement map}.\end{definition}

Although the displacement map $d(\rho_0,\omega_0)=(d_1(\rho_0,\omega_0),d_2(\rho_0,\omega_0))$ is a bidimensional map, there exists a unique analytical function $\tilde{\omega}(\rho_0)$ defined in a neighborhood $V$ of $\rho_0=0$ such that $d_2(\rho_0,\tilde{\omega}(\rho_0))\equiv 0$ (see \cite{GarciaJacobi} for a proof). The function $\mathbf{d}(\rho_0)=d_1(\rho_0,\tilde{\omega}(\rho_0))$, called  \emph{reduced displacement map}, is analytical. Expanding its power series, we have:
$$\mathbf{d}(\rho_0)=v_1\rho_0+v_2\rho_0^2+v_3\rho_0^3+v_4\rho_0^4+\dots.$$
The coefficients $v_k$ are called \emph{focal values} and the coefficients $l_k=v_{2k+1}$ are called \emph{Lyapunov coefficients}. The emphasis in the odd indexed focal values is due to the fact that the first nonzero focal value is the coefficient of an odd power of $\rho_0$ (see \cite{GarciaJacobi}).

The zeros of the reduced displacement map correspond to periodic orbits of system \eqref{eq1}. Moreover, the system has a center on a center manifold $W^c$ if and only if all $l_k$ are null. Thus, the reduced displacement map is a powerful tool to study cyclicity. However, the determination of the focal values is a difficult task.

The study of the cyclicity for three-dimensional systems having a Hopf singularity has a strong relation to the problem of distinguishing if the singular point is either a center or a focus on a center manifold, i.e. the \emph{Center Problem}. The following result provides one of the most useful tools to solve the center problem for system \eqref{eq1}. Its proof can be found in \cite{Bibikov,Adam}.

\begin{theorem}\label{TeoEquiv}
Consider system \eqref{eq1} and let $W^c$ be a center manifold. The following statements are equivalent:
\begin{itemize}
\item[(i)] The origin of the system restricted to $W^c$ is a center;
\item[(ii)] System \eqref{eq1} admits a local analytical first integral $H(x,y,z)$ such that $j^2H(0)=x^2+y^2$;
\item[(iii)] System \eqref{eq1} admits a formal first integral $H(x,y,z)$ such that $j^2H(0)=x^2+y^2$.
\end{itemize}
Moreover, if any of the above statements holds, $W^c$ is unique and analytic.
\end{theorem}

Using Theorem \ref{TeoEquiv}, the standard algorithm to study the center problem is the construction of a formal series $$H(x,y,z)=x^2+y^2+\sum_{j+k+l\geqslant 3}p_{jkl}x^jy^kz^l,$$
with unknown real coefficients $p_{jkl}$. Let $X$ denote the vector field associated to system \eqref{eq1} and consider the following equation:
\begin{equation}\label{XH=0}
XH=\langle X,\nabla H\rangle\equiv 0.
\end{equation}
If it is possible to choose $p_{jkl}$ such that the above equation is satisfied, then \eqref{eq1} has a center on the center manifold (which is, in this case, unique and analytic). Although \eqref{XH=0} is not always satisfied, it is always possible to choose $p_{jkl}$ such that
\begin{equation}\label{XH=Lyap}
XH=\sum_{k\geqslant2}L_{k-1}(x^2+y^2)^k.
\end{equation}
This fact is proved in \cite{Adam}. Furthermore, the quantities $L_{k-1}$ are rational functions whose numerators are polynomials depending on the parameters of system \eqref{eq1}. Any non-zero $L_{k-1}$ is an obstruction for the origin of \eqref{eq1} to be a center on the center manifold. The set of parameters of system \eqref{eq1} for which all $L_{k-1}$ are null, i.e. the origin is a center, is called \emph{Bautin Variety}. The coefficients $L_{k-1}$ are the \emph{Lyapunov constants} for system \eqref{eq1}. The Lyapunov constants and the Lyapunov coefficients are related in the following way: given a positive integer $k>2$, we have
$$l_1=\pi L_{1}\;\mbox{ and }\;l_{k-1}=\pi L_{k-1} \mbox{\hspace{0.3cm} \emph{mod}}\langle L_{1},L_{2},\dots,L_{k-2}\rangle.$$
The proof of this result can be found in \cite{Garcia}. This relationship allows us to obtain the information of the reduced displacement by computing the Lyapunov constants. Since the computations envolving the expressions given in \eqref{XH=Lyap} are algebraic, with the help of symbolic mathematics such as Maple and Mathematica it is possible to compute a large amount the Lyapunov constants with less computational time.

\begin{remark}
The concepts of \emph{Lyapunov constants}, \emph{focal values} and the \emph{Bautin variety} exist for planar systems having an elementary center or focus and their properties and roles are the same as in the three-dimensional case.
\end{remark}

Having the previous tools to obtain the Lyapunov constants, we state the next two results, which can be found in \cite{TorreGouveia} and are proved in \cite{Christopher}.

\begin{theorem}\label{TeoBifLinear}
Suppose that $s$ is a point on the Bautin variety and that the first $k$ Lyapunov coefficients, $L_1,\dots,L_k$, have independent linear parts (with respect to the expansion of $L_i$ about $s$), then $s$ lies on a component of the Bautin variety of codimension at least $k$ and there are bifurcations which produce $k$ limit-cycles locally from the center corresponding to the parameter value $s$. If, furthermore, we know that $s$ lies on a component of the center variety of codimension $k$, then $s$ is a smooth point of the variety, and the cyclicity of the center for the parameter value $s$ is exactly $k$. In the latter case, $k$ is also the cyclicity of a generic point on this component of the Bautin variety.
\end{theorem}

\begin{theorem}\label{TeoBifHighOrder}
Suppose that we are in a point $s$ where Theorem \ref{TeoBifLinear} applies. After a change of variables if necessary, we can assume that $L_0=L_1=\dots=L_k=0$ and the next Lyapunov coefficients $L_i=h_i(u)+O(|u|^{m+1})$, for $i=k+1,\dots,k+l$, where $h_i$ are homogeneous polynomials of degree $m\geqslant 2$ and $u=(u_{k+1},\dots, u_{k+l})$. If there exists a line $\eta$, in the parameter space, such that $h_i(\eta)=0,$ for $i=k+1,\dots,k+l-1$, the hypersurfaces $h_i=0$ intersect transversally along $\eta$ for $i=k+1,\dots, k+l-1$, and $h_{k+l}(\eta)\neq 0$, then there are perturbations of the center which produce $k+l$ limit-cycles.
\end{theorem}

Even though the above theorems where first stated to deal with planar systems, they also apply to the three-dimensional case since they are results regarding the parameter space and the Lyapunov coefficients. 

The approach we use to study the cyclicity for system \eqref{eq1} is as follows: First, we consider systems \eqref{eq1} for which the origin is a center on the center manifold. Then, we make the following perturbation of the considered system:
\begin{equation}\label{eqPerturbation}
\left\{\begin{array}{lcr}
\dot{x}=-y+P(x,y,z)+\sum_{j+k+l=2}a_{jkl}x^jy^kz^l,\\
\dot{y}=x+Q(x,y,z)+\sum_{j+k+l=2}b_{jkl}x^jy^kz^l,\\
\dot{z}=-\lambda z+R(x,y,z)+\sum_{j+k+l=2}c_{jkl}x^jy^kz^l,
\end{array}\right.
\end{equation}
where $\Lambda=(a_{jkl},b_{jkl},c_{jkl})$ for $j+k+l=2$ are the perturbation parameters. Thus, $a_{jkl}=b_{jkl}=c_{jkl}=0$ is a point on the Bautin variety for system \eqref{eqPerturbation}. We then proceed to compute the Lyapunov constants $L_k$ for the perturbed system. Once a sufficient amount $N$ is computed, we evaluate the Jacobian matrix $\frac{\partial (L_1,\dots,L_N)}{\partial (a_{jkl},b_{jkl},c_{jkl})}$ and its rank $r$. We use Theorem \ref{TeoBifLinear} to determine the number of limit-cycles obtainable by the study of the linear part of the Lyapunov constants. 

After this first procedure, we then use Theorem \ref{TeoBifHighOrder}, making the suitable change of variables to verify if the conditions on the hypothesis are satisfied for the quadratic terms of the power series expansion of the Lyapunov constants on the perturbation parameters. If there is a noticeable increase in the number of limit-cycles we repeat this process for the next degree of the expansion until no new limit-cycles are obtained.

\section{R{\"o}ssler, Lorenz and Moon-Rand systems}

\subsection{R{\"o}ssler system:}

In 1976, Rössler proposed the following system
\begin{equation}
	\left\{\begin{array}{lcr}
		\dot{x}=-y-z,\nonumber\\
		\dot{y}=x+az,\nonumber\\
		\dot{z}=b-cz+xz,\nonumber
	\end{array}\right.	
\end{equation}
which has a chaotic behavior for some values of $a,b,c$  \cite{MalasomaRossler}. For parameter values $a=b=0$, the origin is Hopf singular point. The change of variables $x=\bar{x}+\frac{c\bar{z}}{c^2+1}$, $y=\bar{y}-\frac{\bar{z}}{c^2+1}$, $z=\bar{z}$ transforms the above system into
\begin{equation}\label{eqRossler}
	\left\{\begin{array}{lcr}
		\dot{x}=-y-\frac{cxz}{c^2+1}-\frac{c^2z^2}{(c^2+1)^2},\\
		\dot{y}=x+\frac{xz}{c^2+1}+\frac{cz^2}{(c^2+1)^2},\\
		\dot{z}=-cz+xz+\frac{cz^2}{c^2+1},
	\end{array}\right.
\end{equation}
which has a center on the center manifold $W^c=\{z=0\}$ for all real values of $c$. For $c=-1$, we have the following system:
\begin{equation}
	\left\{\begin{array}{lcr}
		\dot{x}=-y+\frac{xz}{2}-\frac{z^2}{4},\nonumber\\
		\dot{y}=x+\frac{xz}{2}+\frac{z^2}{4},\nonumber\\
		\dot{z}=z+xz-\frac{z^2}{2}.\nonumber
	\end{array}\right.
\end{equation}
We compute the first 11 Lyapunov constants. The linear terms of the first $3$ Lyapunov constants are given by
\begin{equation*}
	\begin{aligned}
		&L^1_{1} = -\dfrac{11c_{0, 2, 0}}{15} - \dfrac{c_{1, 1, 0}}{5} - \dfrac{3c_{2, 0, 0}}{5},& \\
		&L^1_{2} = -\dfrac{4c_{0, 2, 0}}{25} - \dfrac{c_{1, 1, 0}}{25} - \dfrac{2c_{2, 0, 0}}{25},& \\
		&L^1_{3} = -\dfrac{101c_{0, 2, 0}}{5950} - \dfrac{3c_{1, 1, 0}}{850} - \dfrac{c_{2, 0, 0}}{170}.& 
	\end{aligned}
\end{equation*}
Their rank is $3$ and thus, by Theorem \ref{TeoBifLinear} it is possible to obtain $3$ limit-cycles. After a suitable coordinate change, we can write $L_i=u_i+O(|\Lambda|^2)$ for $i=1,2,3$ and $L_j=h_j(\Lambda)+O(|\Lambda|^3)$ where $h_j$ is a quadratic polynomial for $j=4,5$. 

\begin{equation*}
	\begin{aligned}
		&h_{4} = \frac{16}{109395}b_{1, 1, 0}b_{2, 0, 0} + \frac{16}{109395}b_{1, 1, 0}b_{0, 2, 0} - \frac{16}{109395}a_{1, 1, 0}a_{2, 0, 0} & \\ 
		&\phantom{=} + \frac{32}{109395}a_{2, 0, 0}b_{2, 0, 0} - \frac{16}{109395}a_{0, 2, 0}a_{1, 1, 0} - \frac{32}{109395}a_{0, 2, 0}b_{0, 2, 0},&  \\
		&h_{5} = \frac{9}{761090}b_{1, 1, 0}b_{2, 0, 0} + \frac{9}{761090}b_{1, 1, 0}b_{0, 2, 0} - \frac{9}{761090}a_{1, 1, 0}a_{2, 0, 0} & \\
		&\phantom{=}+ \frac{9}{380545}a_{2, 0, 0}b_{2, 0, 0} - \frac{9}{761090}a_{0, 2, 0}a_{1, 1, 0} - \frac{9}{380545}a_{0, 2, 0}b_{0, 2, 0}.&
	\end{aligned}
\end{equation*}
Since $h_5$ is a multiple of $h_4$, there is no solution $\Lambda^*$ such that $h_4(\Lambda^*)=0\neq h_5(\Lambda^*)$. Thus, cyclicity of the R{\"o}ssler system \eqref{eqRossler} is at least $4$.

\subsection{Lorenz system:}
The celebrated Lorenz system, first proposed in 1963 \cite{Lorenz}, is one of the widely studied three-dimensional systems due to its rich dynamics. We consider the following generalization of the Lorenz system:
\begin{equation}
\left\{\begin{array}{lcr}
\dot{x}=a(y-x),\nonumber\\
\dot{y}=bx+cy-xz,\nonumber\\
\dot{z}=dz+xz.\nonumber
\end{array}\right.
\end{equation}
For the Lorenz system to have isolated Hopf singularites, we must have $a=c$, $ad\neq 0$ and $a(a+b)<0$. Let $\sigma=\sqrt{-a(a+b)}$. After the suitable coordinate changes are applied, the Lorenz system becomes:
\begin{equation}\label{Lorenz}
	\left\{\begin{array}{lcr}
		\dot{x}=-y-\dfrac{a}{b\sigma}xz+\dfrac{a^2}{b\sigma^2}yz,\\
		\dot{y}=x-\dfrac{1}{b}xz+\dfrac{a}{b\sigma}yz,\\
		\dot{z}=\dfrac{d}{\sigma}z+\dfrac{1}{b}xy-\dfrac{a}{b\sigma}y^2.
	\end{array}\right.
\end{equation}

The Bautin variety for system \eqref{Lorenz} is given by the condition $d=-2a$. Under those, considering the perturbation \eqref{eqPerturbation} we compute the Lyapunov constants for system \eqref{Lorenz}.

For parameter values $a=-1,b=5,d=2$, system \eqref{Lorenz} becomes:
\begin{equation}
\left\{\begin{array}{lcr}
\dot{x}=-y+\dfrac{1}{10}xz+\dfrac{1}{20}yz,\nonumber\\
\dot{y}=x-\dfrac{1}{5}xz-\dfrac{1}{10}yz,\nonumber\\
\dot{z}=z+\dfrac{1}{5}xy+\dfrac{1}{10}y^2.\nonumber
\end{array}\right.
\end{equation}
The linear part of the first eleven Lyapunov constants has rank $2$ and by Theorem \ref{TeoBifLinear} it is possible to obtain $2$ limit-cycles from the origin. Making an appropriate coordinate change, we can write $L_1=u_1+O(|\Lambda|^2)$, $L_2=u_2+O(|\Lambda|^2)$ and $L_3=h_3(\Lambda)+O(|\Lambda|^3)$ and $L_4=h_4(\Lambda)+O(|\Lambda|^3)$ where $h_3,h_4$ are quadratic polynomials. It is possible to find a solution $\Lambda^*$ such that $L_3(\Lambda^*)=0\neq L_4(\Lambda^*)$. Using Theorem \ref{TeoBifHighOrder}, it is possible to obtain two additional limit-cycles from the center at the origin. We conclude that the cyclicity of the Lorenz system is at least $4$. 

\subsection{Moon-Rand system:}

The Moon-Rand \cite{MoonRand} system is given by
\begin{equation}\label{Moon-Rand}
\left\{\begin{array}{lcr}
\dot{x}=y,\\
\dot{y}=-x-xz,\\
\dot{z}=-\mu z+cx^2+bxy+ay^2,
\end{array}\right.
\end{equation}
and its Bautin variety is $V^c={a=2c-\mu b=0}$, see \cite{Queiroz,Garcia}. Assuming the follow parameter values $\mu=1,b=2,c=1$, we obtain the next system.
\begin{equation}
\left\{\begin{array}{lcr}
\dot{x}=y,\nonumber\\
\dot{y}=-x-xz,\nonumber\\
\dot{z}=-z+x^2+2xy,\nonumber
\end{array}\right.	
\end{equation}
The linear part of the first ten Lyapunov constants has rank $2$ and by Theorem \ref{TeoBifLinear} it is possible to obtain $2$ limit-cycles from the origin. Making an appropriate coordinate change, we can write $L_1=u_1+O(|\Lambda|^2)$, $L_2=u_2+O(|\Lambda|^2)$ and $L_3=h_3(\Lambda)+O(|\Lambda|^3)$ and $L_4=h_4(\Lambda)+O(|\Lambda|^3)$ where $h_3,h_4$ are quadratic polynomials. As in the Lorenz system, it is possible to find a solution $\Lambda^*$ such that $L_3(\Lambda^*)=0\neq L_4(\Lambda^*)$. By Theorem \ref{TeoBifHighOrder}, it is possible to obtain two additional limit-cycles from the center at the origin. We conclude that the cyclicity of the Moon-Rand system is at least $4$.

The results of this section are summed up in Theorem \ref{TeoLorenzRosslerMoonRand}.		
 
\section{Jerk systems}

Consider the three-dimensional system of differential equations:
\begin{equation}\label{eqJerk}
\left\{\begin{array}{lcr}
\dot{x}=y,\\
\dot{y}=z,\\
\dot{z}=f(x,y,z).
\end{array}
\right.
\end{equation}
It is equivalent to the third-order differential equation $\dddot{x}=f(x,\dot{x},\ddot{x})$ which is called \emph{Jerk equation} since in mechanical models where $x(t)$ denotes the position of given particle at instant $t$, $\dddot{x}(t)$ denotes the rate of change of the acceleration of such particle, that is, its jerk. Jerk equations are broadly studied in the literature for its use in science and engineering. 

We search for Hopf singular points for system \eqref{eqJerk} with $f(x,y,z)$ being a polynomial. The singular points, when they exist, are given by $Q_i=(\alpha_i,0,0)$ where $\alpha_i$ are the roots of $f(x,0,0)$ for $i=1,\dots,\deg(f)$. The Jacobian matrix $J$ of \eqref{eqJerk} is always given by:
$$J=\left(\begin{array}{lcr}
	0 & 1 & 0 \\
	0 & 0 & 1 \\
	\dfrac{\partial f}{\partial x} & \dfrac{\partial f}{\partial y} &\dfrac{\partial f}{\partial z} 
\end{array}\right).$$
The determinant and the trace of the above matrix are, respectively, $\det J=\frac{\partial f}{\partial x}$ and ${\rm tr} J=\frac{\partial f}{\partial z}$. Its characteristic polynomial is the following:
$$-\eta^3+\eta^2\frac{\partial f}{\partial z}+\eta\frac{\partial f}{\partial y}+\frac{\partial f}{\partial x}.$$
For $Q_i$ to be a nilpotent singular point, we must have $\frac{\partial f}{\partial z}(Q_i)=\tau$, $\frac{\partial f}{\partial y}(Q_i)=-\beta^2$ and $\frac{\partial f}{\partial x}(Q_i)=\beta^2\tau$ for some $\beta,\tau\neq 0$.

Translating $Q_{i}$ to the origin transforms system \eqref{eqJerk} into:
\begin{equation}
\left\{\begin{array}{lcr}
\dot{x}=y,\nonumber\\
\dot{y}=z,\nonumber\\
\dot{z}=f(x+\alpha_i,y,z).
\end{array}\right.	
\end{equation}
Performing the change of variables $x=-\frac{\bar{x}}{\beta^2}+\frac{\bar{z}}{\tau^2},\;y=\frac{\bar{y}}{\beta}+\frac{\bar{z}}{\tau}$, $z=\bar{x}+\bar{z}$, dropping the bars, and rescaling time, we obtain the canonical form:

\begin{equation}\label{eqJerk2}
	\left\{\begin{array}{lcr}
	\dot{x}=- y+\beta F_2(x,y,z),\\
	\dot{y}=x-\tau F_2(x,y,z),\\
	\dot{z}=-\lambda z-\tau\lambda F_2(x,y,z),	
	\end{array}\right.	
\end{equation}
where $F_2(x,y,z)$ has no linear nor constant terms. We further investigate the origin of system \eqref{eqJerk2} since it is the canonical form of any jerk system having a Hopf singularity. An extensive study on the center problem for some families of Jerk systems can be found in \cite{MahdiJerk, MahdiPessoa}. We study the cyclicity of the centers given by each of the center conditions found in these papers.

We consider $F_2$ a quadratic polynomial, i.e. $F_2(x,y,z)=a_1x^2+a_2y^2+a_3z^2+a_4xy+a_5xz+a_6yz$ and $\lambda=\beta=-\tau=1$. The center conditions for those systems, proven by the authors of \cite{MahdiPessoa}, are:
\begin{itemize}
	\item[a)]$a_1=a_2=a_4=0$;
	\item[b)]$a_1-a_2=a_3=a_5=a_6=0$;
	\item[c)]$a_1+a_2=a_3=a_5=a_6=0$;
	\item[d)]$a_1+a_2=2a_2-a_3+a_6=a_3-a_4-2a_5=2a_4+3a_5+a_6=0$;
	\item[e)]$2a_1-a_6=2a_2+a_5=2a_3-a_5+a_6=a_4+a_5+a_6=0$;
	\item[f)]$a_1-a_2=2a_2+a_6=a_4=a_5+a_6=0$;
	\item[g)]$2a_1+a_2=2a_2+a_6=4a_3+5a_6=a_4=2a_5-a_6=0$.
\end{itemize}

For the case $\textbf{a)}$, we take the parameter values $a_3=1,a_5=2,a_6=3$. Then, perturbing the system with 18 perturbation parameters, we obtain 3 limit-cycles using the linear terms of the Lyapunov constants. In the case $\textbf{b)}$, we assume $a_2=2, a_4=4$ and obtain 5 limit-cycles with the linear terms of Lyapunov constants.

Analyzing the case  $\textbf{c)}$ and assuming again $a_2=2, a_4=4$, we obtain 6 limit-cycles using linear terms. In the case $\textbf{d)}$, taking $a_4=1,a _5=-1/2$, we also obtain 5 limit-cycles using only the linear terms of Lyapunov constants.

For case $\textbf{e)}$, assuming $a_2=5,a_3=-1/2$, computing the linear terms of the first 11 Lyapunov constants, we obtain 6 limit-cycles.

For the case $\textbf{f)}$, taking $a_2=3/5,a_3=1$, we obtain the system 
\begin{equation*}\left\{\begin{array}{lll}
\dot{x} = -y+\frac{3}{5}x^2 + \frac{6}{5}xz + \frac{3}{5}y^2 - \frac{6}{5}yz + z^2, \vspace{0.25cm} \\ 
\dot{y} = x+\frac{3}{5}x^2 + \frac{6}{5}xz + \frac{3}{5}y^2 - \frac{6}{5}yz + z^2, \vspace{0.25cm} \\
\dot{z} =-z+\frac{3}{5}x^2 + \frac{6}{5}xz + \frac{3}{5}y^2 - \frac{6}{5}yz + z^2.
\end{array}
\right.
\end{equation*}

Computing the first eleven Lyapunov constants up to order 2 for this system, we obtain 7 limit-cycles using the linear terms with the perturbation parameters $a_{0, 0, 2}$, $a_{0, 1, 1}$, $a_{0, 2, 0}$, $a_{1, 0, 1}$, $a_{1, 1, 0}$, $a_{2, 0, 0}$, $b_{0, 1, 1}$. Making an appropriate change of coordinates to vanish the linear terms of the Lyapunov constants $L_{8}$ through $L_{11}$, we obtain one more limit cycle. Therefore, it is possible to obtain at least 8 limit-cycles with Lyapunov constants up to order 2. 

\medskip

\noindent\textbf{Proof of Theorem \ref{TeoJerk12}: }
Finally, among the considered jerk systems the one with most cyclicity was system \eqref{eqJerk12}, which we recall here: 
\begin{equation*}\left\{\begin{array}{lll}
	\dot{x}&=&-y+2x^2+4xz-4y^2+8yz-10z^2,\vspace{0.25cm} \\ 
	\dot{y}&=&x+2x^2+4xz-4y^2+8yz-10z^2,\vspace{0.25cm} \\ 
	\dot{z}&=&-z+2x^2+4xz-4y^2+8yz-10z^2,
\end{array}
	\right.
\end{equation*}
The origin is a center on the center manifold since the system satisfies center condition \textbf{g)} above. To simplify the calculations, we assume $b_{0, 2, 0}= b_{1, 0, 1}= b_{1, 1, 0}= b_{2, 0, 0}= c_{0, 0, 2}=c_{0, 1, 1}=0$ in the quadratic perturbation \eqref{eqPerturbation}.  The rank of the linear part of the first 12 Lyapunov constants is $8$ with the perturbation parameters $a_{0, 0, 2}, a_{0, 1, 1}, a_{0, 2, 0}, a_{1, 0, 1}$, $a_{1, 1, 0}, a_{2, 0, 0}, b_{0, 0, 2}, b_{0, 1, 1}$. 

We make following change of variables: $a_{0, 0, 2}=\gamma^{2}u_{1}$, $a_{0, 1, 1}=\gamma^{2}u_{2}$, $ a_{0, 2, 0}=\gamma^{2}u_{3}$, $a_{1, 0, 1}=\gamma^{2}u_{4}$, $a_{1, 1, 0}=\gamma^{2}u_{5}$, $a_{2, 0, 0}=\gamma^{2}u_{6}$, $ b_{0, 0, 2}=\gamma^{2}u_{7}$, $b_{0, 1, 1}=\gamma^{2}u_{8}$, $c_{0, 2, 0} = u_{9}\gamma$, $c_{1, 0, 1} = u_{10}\gamma$, $c_{1, 1, 0} = u_{11}\gamma$, $c_{2, 0, 0} = \gamma$. Considering the 2nd order terms Taylor expansion of the Lyapunov constants in $\gamma$, we obtain eleven equations in eleven variables. Due to the size of equations, we show only the first three of them:

\begin{align*}
&\mathcal{L}_{1} = \dfrac{107773}{25992}u_{10} - \dfrac{37355}{12996}u_{11}^2 - \dfrac{3869}{4332}u_{10}^2 - \dfrac{3283}{25992}u_{9}^2  + \dfrac{30829}{2888}u_{11} + \dfrac{7377}{1444}u_{9} \\ & \phantom{=}+ u_{1} + \dfrac{853}{25992} - \dfrac{110623}{25992}u_{10}u_{11} - \dfrac{11047}{25992}u_{9}u_{10} + \dfrac{385}{25992}u_{9}u_{11}, \\
&\mathcal{L}_{2} = \dfrac{10472438}{81225}u_{11}^2 - \dfrac{1405238}{81225}u_{10}^2 - \dfrac{203501}{16245}u_{9}^2 + \dfrac{2717813}{27075}u_{10}u_{11} + u_{2} \\ & \phantom{=}- \dfrac{72276202}{81225}u_{9} + \dfrac{13572899}{81225}u_{11} + \dfrac{5841047}{16245}u_{10} + \dfrac{3562733}{27075}u_{9}u_{10}  + \dfrac{728603}{16245}u_{9}u_{11}\\ & \phantom{=} - \dfrac{10443633}{9025}, \\
&\mathcal{L}_{3}=-\dfrac{64427829146}{9665775}u_{10}u_{11} + u_{3} + \dfrac{35151434924}{386631}u_{9} - \dfrac{1119816144458}{48328875}u_{11} \\ & \phantom{=} - \dfrac{218875710314}{5369875}u_{10} - \dfrac{167865339164}{16109625}u_{11}^2 + \dfrac{27891420988}{9665775}u_{10}^2 + \dfrac{6274007402}{3221925}u_{9}^2 \\ & \phantom{=} - \dfrac{693744306994}{48328875}u_{9}u_{10} - \dfrac{37253310546}{5369875}u_{9}u_{11} + \dfrac{1168498101998}{9665775}.
\end{align*}
 
Solving $\mathcal{L}_{1},\cdots,\mathcal{L}_{11}$ in the variables $u_{i},i=1,\cdots,11$, we obtain 2 solutions. The first solution give us a center. The second solution is given by

\begin{align*}
&u_{1}= -\dfrac{85781\cdots88188}{99469\cdots13561}, & u_{2} =\dfrac{43811\cdots03456}{31538\cdots48581},\\ & 
u_{3}=-\dfrac{26274\cdots38332}{31538\cdots48581}, & u_{4}=\dfrac{54792\cdots 91712}{24867\cdots39025}, \\ & 
u_{5} = -\dfrac{34934\cdots 50976}{14796\cdots19875}, & u_{6} =\dfrac{12648\cdots 64608}{32699\cdots92375}, \\ &
u_{7}= -\dfrac{66960\cdots29504}{73530\cdots53125}, & u_{8} =\dfrac{11188\cdots98112}{39081\cdots03125}, \\ & 
u_{9} =-\dfrac{41701\cdots20512}{35951\cdots65625}, & u_{10}=\dfrac{30214\cdots72768}{51225\cdots65625}, \\ &
u_{11} = \dfrac{80438\cdots 84083}{31538\cdots48581}.
\end{align*}
Evaluating $\mathcal{L}_{1},\cdots,\mathcal{L}_{12}$ at this solution gives us $\mathcal{L}_{1}=0,\cdots,\mathcal{L}_{11}=0$ and
 \[\mathcal{L}_{12}=\dfrac{12619\cdots11136}{70884\cdots59375}.\]
 Moreover, the Jacobian matrix $J$ of $\mathcal{L}_{1},\cdots,\mathcal{L}_{11}$ with respect to $u_{i}$, $i=1,\cdots,11$ at the above solution, has the following determinant 
 \[\det J=\dfrac{52981\cdots61344}{18459\cdots90625}.\] 

Therefore, by Theorem \ref{TeoBifHighOrder}, 12 limit-cycles unfold from the center. This proves Theorem \ref{TeoJerk12}.

\section{Gin\'e-Valls systems}

The authors of \cite{Gine} found conditions on the parameters for which the origin of the system \eqref{Valls}, i.e.
\begin{equation*}
	\left\{\begin{array}{lcr}
	\dot{x}=y,\\
	\dot{y}=-x+a_1x^2+a_2xy+a_3xz+a_4y^2+a_5yz+a_6z^2,\\
	\dot{z}=-z+c_1x^2+c_2xy+c_3y^2,	
	\end{array}\right.	
\end{equation*}
is a center. The conditions found in the cited paper are in the Bautin Variety of \eqref{Valls}. We go through each condition and implement our study of the cyclicity, computing the Lyapunov constants for a perturbation \eqref{eqPerturbation} for system \eqref{Valls}.

\subsection{Case $a_1=a_2=0$:}

For system \eqref{Valls} with $a_1=a_2=0$, the origin is a center if one of the following conditions holds:
\begin{itemize}
	\item[a)]$c_1=c_2=c_3=0$;
	\item[b)]$a_5=2c_1-c_2=c_3=0$;
	\item[c)]$a_3=a_5=a_6=0$.
\end{itemize}

In items \textbf{a)} and \textbf{b)}, assuming $a_5=3/4,a_3=1,a_4=3,a_6=1/2$ and $c_1=0,a_3=4,a_4=1/2,a_6=-5/3$ respectively,  we obtain 4 limit-cycles with the linear terms of Lyapunov constants for both systems. For the item \textbf{c)}, taking $c_1=-1,c_2=2/3,c_3=5,a_4=-2/3$ we obtain 7 limit-cycles.

\subsection{Case $a_1=a_3=0$:}

For system \eqref{Valls} with $a_1=a_3=0$, the origin is a center if one of the following conditions holds:
\begin{itemize}
	\item[a)]$a_4=c_1=c_2=c_3=0$;
	\item[b)]$a_2=c_1=c_2=c_3=0$;
	\item[c)]$2a_{2}a_{4}+a_{5}c_{2}=a_{4}a_{5}^{2}-a_{2}^2a_{6}=2a_{4}^2a_{5}+a_{2}a_{6}c_2=4a_{4}^3-a_{6}c_{2}^2=2a_{2}^3 a_{6}+a_{5}^3 c_{2}=2c_1-c2=c_3=0$;
	\item[d)]$a_2=a_5=2c_1-c_2=c_3=0$;
	\item[e)]$a_4=a_5=a_6=0$;
	\item[f)]$a_2=a_5=a_6=0$;
\end{itemize}

In the item \textbf{a)}, considering $a_2=3,a_5=2/3,a_6=-6$, we obtain 3 limit-cycles. For the item \textbf{b)}, assuming $a_4=5,a_5=2,a_6=1$ we obtain 4 limit-cycles.

Analyzing the item \textbf{c)}, taking $a_2=-1,a_4=5/8,c_1=1/2$, we obtain the system 

\begin{equation*}\left\{\begin{array}{lll}
\dot{x} = y, \vspace{0.25cm} \\ 
\dot{y} = -xy + \dfrac{5}{8}y^2 + \dfrac{5}{4}yz + \dfrac{125}{128}z^2-x, \vspace{0.25cm} \\
\dot{z} =\dfrac{1}{2}x^2 + xy - z.
\end{array}
\right.
\end{equation*}
Computing the Lyapunov constants up to order 3, we obtain using the linear terms 9 limit-cycles with the perturbation parameters $a_{0,0,2},a_{0,1,1}$, $a_{0,2,0},a_{1,0,1},a_{1,1,0},a_{2,0,0},b_{0,0,2},b_{0,1,1},c_{0, 2, 0}$. Then, via a change of variables, we can write the first nine Lyapunov constants as follows

\begin{equation*}\left\{\begin{array}{lll}
L_{1} & = & u_{1} + O(|\Lambda|^2), \\ 
L_{2} & = & u_{2} + O(|\Lambda|^2),  \\
\vdots&   &  \\
L_{9} & = & u_{9} + O(|\Lambda|^2). \\
\end{array}
\right.
\end{equation*}
Doing an appropriate change of coordinates and vanishing the linear terms of $L_{10}$ and $L_{11}$ we have

\begin{equation*}\left\{\begin{array}{lll}
L_{10} & = & u_{10}u_{11} + O(|\Lambda|^3), \\ 
L_{11} & = & u_{10}u_{11} + O(|\Lambda|^3),  \\
\end{array}
\right.
\end{equation*}
where $u_{10}=c_{0, 0, 2}$ and $u_{11}=8c_{0, 0, 2} + 4c_{1, 0, 1} - 9c_{0, 1, 1} + 5c_{1, 1, 0} + 4b_{1, 0, 1} - 10c_{2, 0, 0} - 10b_{2, 0, 0}$. So, clearly, we can obtain one more limit-cycle. Therefore, when perturbing the system, we can obtain 10 limit-cycles. No more limit-cycles are obtained with Lyapunov constants up to order 3.

In the item \textbf{d)}, taking $a_4=1,a_6=-1,c_1=2/3$, we obtain de system 

\begin{equation*}\left\{\begin{array}{lll}
\dot{x} = y, \vspace{0.25cm} \\ 
\dot{y} = y^2 - z^2 - x, \vspace{0.25cm} \\
\dot{z} =\dfrac{2}{3}x^2 + \dfrac{4}{3}xy - z.
\end{array}
\right.
\end{equation*}
Computing the Lyapunov constants up to order 3, we have 8 limit-cycles with the linear terms, 1 limit cycle from the terms of order 2 and no more limit-cycles with third order terms. We can write the Lyapunov constants as
\begin{equation*}\left\{\begin{array}{lll}
L_{1} & = & u_{1} + O(|\Lambda|^2), \\ 
L_{2} & = & u_{2} + O(|\Lambda|^2),  \\
\vdots&   &  \\
L_{9} & = & u_{9}u_{10} + O(|\Lambda|^3), \\
L_{9} & = & u_{9}u_{10} + O(|\Lambda|^3), \\
\end{array}
\right.
\end{equation*}
where $u_{9}=c_{0, 0, 2}$ and $u_{10}=2c_{1, 0, 1} + 3c_{0, 2, 0} - 2c_{0, 1, 1} + 2b_{1, 1, 0}$.

\subsection{Case $a_1=a_4=0$:}

For system \eqref{Valls} with $a_1=a_4=0$, the origin is a center if one of the following conditions holds:
\begin{itemize}
\item[a)]$c_1=c_2=c_3=0$;
\item[b)]$a_5=a_6=2c_1-c_2=c_3=0$;
\item[c)]$a_2=a_5=2c_1-c_2=c_3=0$;
\item[d)]$a_3=a_5=a_6=0$.
\end{itemize}
For the item \textbf{a)}, taking parameter values $a_{2}=1, a_{3}=2, a_{5}=-1/2$ and $a_{6}=5/3$ we obtain 3 limit-cycles using the linear terms of Lyapunov constants. In item \textbf{b)}, for $c_{1}=1/2,
a_{2}=-2$ and $a_{3}=7/8$ we obtain 5 limit-cycles with linear terms. 

Studying item \textbf{c)}, this condition produced most limit-cycles for this case. Taking $a_{3}=7/8,
a_{6}=2/8$ and $c_{1}=2/3$ we obtain the system

\begin{equation*}\left\{\begin{array}{lll}
\dot{x} = y, \vspace{0.25cm} \\ 
\dot{y} = \dfrac{7}{8}xz+\dfrac{1}{4}z^2-x, \vspace{0.25cm} \\
\dot{z} =\dfrac{2}{3}x^2+\dfrac{4}{3}xy-z.
\end{array}
\right.
\end{equation*}
Computing the first 11 Lyapunov constants of order 3, we obtain rank $7$. Then, by the Theorem \ref{TeoBifLinear}, we have, at least,  $7$ limit-cycles with the perturbation parameters $a_{0,0,2},a_{0,1,1},a_{0,2,0},a_{1,0,1},a_{1,1,0},a_{2,0,0}$,$c_{0, 2, 0}$. Doing an appropriate change of coordinates, we can write the Lyapunov constants as 
\begin{equation*}\left\{\begin{array}{lll}
L_{1} & = & u_{1} + O(|\Lambda|^2), \\ 
L_{2} & = & u_{2} + O(|\Lambda|^2),  \\
\vdots&   &  \\
L_{7} & = & u_{7} + O(|\Lambda|^2), \\
L_{8} & = & u_{8}u_{9}+ O(|\Lambda|^3), \\
L_{9} & = & u_{8}u_{10} + O(|\Lambda|^3), \\
L_{10} & = & u_{8}u_{9} + O(|\Lambda|^3), \\
\end{array}
\right.
\end{equation*}
where $u_{8}=32b_{0, 1, 1} - 7c_{1, 1, 0} + 14c_{2, 0, 0}$, $u_{9}$ and $u_{10}$ are homogeneous polynomial of degree 2 in the variables $b_{0, 1, 1}$, $b_{0, 2, 0}$, $b_{1, 1, 0}$, $c_{0, 0, 2}$, $c_{0, 1, 1}$, $c_{1, 0, 1}$, $c_{1, 1, 0}$, $c_{2, 0, 0}$. So, it is clear that we can obtain 9 limit-cycles.

Finally, for the item \textbf{d)}, taking $a_{2}=1/2,c_{1}=-2,c_{2}=7/8,c_{3}=2/8$ we obtain 6 limit-cycles using the linear terms of the Lyapunov constants.

\subsection{Case $a_1=a_5=0$:}
 
For system \eqref{Valls} with $a_1=a_5=0$, the origin is a center if one of the following conditions holds:
\begin{itemize}
 	\item[a)]$a_4=c_1=c_2=c_3=0$;
 	\item[b)]$a_4=a_6=2c_1-c_2=c_3=0$;
 	\item[c)]$a_2=2c_1-c_2=c_3=0$;
 	\item[d)]$a_3=a_4=a_6=0$;
 	\item[e)]$a_2=a_3=a_6=0$.
\end{itemize}
In the item \textbf{a)}, taking $a_{2}=1/2 a_{3}=-2,a_{6}=7/8$, we obtain 3 limit-cycles with the linear terms of Lyapunov constants. For the item  \textbf{b)}, we have 5 limit-cycles taking $a_{3}=2,a_{2}=3,c_{1}=-2/3$.

Analyzing the item \textbf{c)} with $a_{3}=-2,a_{4}=-2,a_{6}=7/8,c_{1}=2/3$ we have the system 
 
 \begin{equation*}\left\{\begin{array}{lll}
 \dot{x} = y, \vspace{0.25cm} \\ 
 \dot{y} = -2xz - 2y^2 + \dfrac{7}{8}z^2 - x, \vspace{0.25cm} \\
 \dot{z}=\dfrac{2}{3}x^2 + \dfrac{4}{3}xy-z.
 \end{array}
 \right.
 \end{equation*}
 
 Computing the first 11 Lyapunov constants of order 2, we obtain 8 limit-cycles using linear terms of Lyapunov constants with the parameters  $a_{0, 0, 2}, a_{0, 1, 1}, a_{0, 2, 0}, a_{1, 0, 1}, a_{1, 1, 0}, a_{2, 0, 0}, b_{0, 1, 1},c_{2, 0, 0}$. After an appropriate change of coordinates, we can write the Lyapunov constants as 
 
 \begin{equation*}\left\{\begin{array}{lll}
 L_{1} & = & u_{1} + O(|\Lambda|^2), \\ 
 L_{2} & = & u_{2} + O(|\Lambda|^2),  \\
 \vdots&   &  \\
 L_{8} & = & u_{8} + O(|\Lambda|^2), \\
 L_{9} & = & u_{9}u_{10}+ O(|\Lambda|^3), \\
 L_{10} & = & u_{9}u_{10} + O(|\Lambda|^3), \\
 \end{array}
 \right.
 \end{equation*}
where $c_{0, 0, 2}=u_{9}$ and $-c_{0, 1, 1} - 3c_{0, 2, 0} + c_{1, 0, 1} + b_{1, 1, 0}=u_{10}$. Therefore, we can obtain 9 limit-cycles, being 8 limit-cycles of linear terms of Lyapunov constants and 1 limit cycle with the 2nd order terms of the Lyapunov constants.

For item \textbf{d)}, taking $a_{2}=1/2, c_{3}=2/8, c_{1}=2/3,c_{2}=-1/2$, we obtain 6 limit-cycles using the linear terms of Lyapunov constants. In the item \textbf{e)}, we have 7 limit-cycles assuming $a_{4}=-2, c_{3}=2/8,c_{1}=2/3, c_{2}=-1/2$.

\subsection{Case $a_1=a_6=0$:}

For system \eqref{Valls} with $a_1=a_6=0$, the origin is a center if one of the following conditions holds:
\begin{itemize}
	\item[a)]$a_4=c_1=c_2=c_3=0$;
	\item[b)]$a_2=c_1=c_2=c_3=0$;
	\item[c)]$a_4=a_5=2c_1-c_2=c_3=0$;
	\item[d)]$a_2=a_5=2c_1-c_2=c_3=0$;
	\item[e)]$a_3=a_4=a_5=0$;
	\item[f)]$a_2=a_3=a_5=0$;
\end{itemize}

Analyzing item \textbf{a)}, if we take $a_{2}=1/2,a_{3}=-2,a_{5}=7/8$, we obtain 3 limit-cycles with the linear terms os Lyapunov constants. For item  \textbf{b)}, assuming $a_{3}=-2, a_{4}=-2, a_{5}=7/8$, we obtain 4 limit-cycles. In the idem \textbf{c)}, taking $a_{2}=1/2,a_{3}=-2,c_{1}=2/3$ we have 5 limit-cycles.

For item \textbf{d)}, assuming $a_{3}=-2, a_{4}=-2,c_{1}=2/3$, we obtain the system 
\begin{equation*}\left\{\begin{array}{lll}
\dot{x} = y, \vspace{0.25cm} \\ 
\dot{y} = -2xz-2y^2-x, \vspace{0.25cm} \\
\dot{z} =\dfrac{2}{3}x^2+\dfrac{4}{3}xy-z.
\end{array}
\right.
\end{equation*}
This case follow the same steps from the item \textbf{c)} for the case $a_{1}=a_{5}=0$. Then, we have 

\begin{equation*}\left\{\begin{array}{lll}
L_{1} & = & u_{1} + O(|\Lambda|^2), \\ 
L_{2} & = & u_{2} + O(|\Lambda|^2),  \\
\vdots&   &  \\
L_{8} & = & u_{8} + O(|\Lambda|^2), \\
L_{9} & = & u_{9}u_{10}+ O(|\Lambda|^3), \\
L_{10} & = & u_{9}u_{10} + O(|\Lambda|^3), \\
\end{array}
\right.
\end{equation*}
where, $u_{9}=c_{0, 0, 2}$ and $u_{10}=c_{1, 0, 1} - 3c_{0, 2, 0} - c_{0, 1, 1} + b_{1, 1, 0}$. Therefore, we have 9 limit-cycles bifurcating from the center.

For the items \textbf{e)} and \textbf{f)}, taking $a_{3}=0,a_{4}=0,a_{5}=0$ and $a_{2}=0,a_{3}=0,a_{5}=0$ respectively we obtain 6 and 7 limit-cycles respectively.

\subsection{Case $a_2=a_3=0$:}

For system \eqref{Valls} with $a_2=a_3=0$, the origin is a center if one of the following conditions holds:
\begin{itemize}
	\item[a)]$a_5=a_6=0$;
	\item[b)]$a_5=2c_1-c_2=c_3=0$;
	\item[c)]$c_1=c_2=c_3=0$.
\end{itemize}

In the item \textbf{a)}, taking $a_{1}=1/2,a_{4}=-2,c_{3}=2/8,c_{1}=2/3,c_{2}=-1/2$, we obtain the system 

\begin{equation*}\left\{\begin{array}{lll}
\dot{x} = y, \vspace{0.25cm} \\ 
\dot{y} = \dfrac{1}{2}x^2-2y^2-x, \vspace{0.25cm} \\
\dot{z} =\dfrac{2}{3}x^2 - \dfrac{1}{2}xy + \dfrac{1}{4}y^2-z.
\end{array}
\right.
\end{equation*}

Computing the Lyapunov constants up to order 3, it follows that with linear terms, we have rank 8 with the perturbation parameters $a_{0, 0, 2}$, $a_{0, 1, 1}$, $a_{0, 2, 0}$, $a_{1, 0, 1}$, $a_{2, 0, 0}$, $b_{0, 0, 2}$, $b_{0, 1, 1}$, $b_{1, 0, 1}$. Studying terms of order 2 and 3, we obtain no more limit-cycles.

This case follow the same steps from the item \textbf{c)} for the case $a_{1}=a_{5}=0$. Then, we have 

\begin{equation*}\left\{\begin{array}{lll}
L_{1} & = & u_{1} + O(|\Lambda|^2), \\ 
L_{2} & = & u_{2} + O(|\Lambda|^2),  \\
\vdots&   &  \\
L_{8} & = & u_{8} + O(|\Lambda|^2), \\
L_{9} & = & u_{9}u_{10}+ O(|\Lambda|^3), \\
L_{10} & = & u_{9}u_{10} + O(|\Lambda|^3), \\
\end{array}
\right.
\end{equation*}
where, $u_{9}=c_{0, 0, 2}$ and $u_{10}=-4c_{0, 1, 1} - 15c_{0, 2, 0} + 4c_{1, 0, 1} + 4b_{1, 1, 0}$. Therefore, we have 9 limit-cycles bifurcating from the center, being 8 with linear terms of Lyapunov constants and 1 with of terms of order 2.

For item \textbf{c)}, assuming $a_{1}=1/2,a_{4}=-2,a_{5}=7/8,a_{6}=1$, we obtain 4 limit-cycles.

\subsection{Case $a_2=a_4=0$:}

For system \eqref{Valls} with $a_2=a_4=0$, the origin is a center if one of the following conditions holds:
\begin{itemize}
	\item[a)]$a_5=2c_1-c_2=c_3=0$;
	\item[b)]$c_1=c_2=c_3=0$; 
	\item[c)]$a_3=a_5=a_6=0$; 
	\end{itemize}

For this case, we have for each item 7, 3 and 5 limit-cycles respectively. To obtain this values, we take $a_1=1/2,a_3=-2,a_6=1,c_1=2/3$, $a_1=1/2,a_3=-2,a_5=7/8,a_6=1$ and $a_1=1/2,c_3=2/8,c_1=2/3,c_2=-1/2$ respectively.

\subsection{Case $a_2=a_5=0$:}

For system \eqref{Valls} with $a_2=a_5=0$, the origin is a center if one of the following conditions holds:
\begin{itemize}
	\item[a)]$a_3=a_6=0$;
	\item[b)]$a_1=a_3=a_4=c_1=c_2=c_3=0$; 
	\item[c)]$2c_1-c_2=c_3=0$; 
\end{itemize}

In the item \textbf{a)}, assuming $a_1=1/2,a_4=-2,c_3=2/8,c_1=2/3,c_2=-1/2$, we obtain the system
\begin{equation*}\left\{\begin{array}{lll}
\dot{x} = y, \vspace{0.25cm} \\ 
\dot{y} = \dfrac{1}{2}x^2 - 2y^2 - x, \vspace{0.25cm} \\
\dot{z} =\dfrac{2}{3}x^2 - \dfrac{1}{2}xy + \dfrac{1}{4}y^2 - z.
\end{array}
\right.
\end{equation*}
When perturbirng this system, with the linear terms of Lyapunov constants, we obtain 8 limit-cycles using the perturbation parameters $a_{0, 0, 2}, a_{0, 1, 1}, a_{0, 2, 0}, a_{1, 0, 1}, a_{2, 0, 0}, b_{0, 0, 2}, b_{0, 1, 1},b_{1, 0, 1}$. After an appropriate change of variables, we obtain no more limit-cycles with terms of order 2 and 3 of the Lyapunov constants.

The item \textbf{b)} taking $a_6=1$, we compute the linear terms of the first 10 Lyapunov constants and we obtain rank 0.

Analyzing item \textbf{c)} and taking $a_1=1/2,a_3=-2,a_4=-2,a_5=7/8,a_6=1,c_3=2/8,c_1=2/3,c_2=-1/2$, we obtain the system

\begin{equation*}\left\{\begin{array}{lll}
\dot{x} = y, \vspace{0.25cm} \\ 
\dot{y} = \dfrac{1}{2}x^2 - 2xz - 2y^2 + z^2 - x, \vspace{0.25cm} \\
\dot{z} =\dfrac{2}{3}x^2 + \dfrac{4}{3}xy - z.
\end{array}
\right.
\end{equation*}
With the perturbation parameters
$a_{0, 0, 2}$, $a_{0, 1, 1}$, $a_{0, 2, 0}$, $a_{1, 0, 1}$, $a_{1, 1, 0}$, $a_{2, 0, 0}$, $b_{0, 1, 1}$, $c_{2, 0, 0}$, using the linear terms of the Lyapunov constants, we obtain 8 limit-cycles. Then, doing an appropriate chance of coordinates, we can write the Lyapunov constants as

\begin{equation*}\left\{\begin{array}{lll}
L_{1} & = &u_{1} + O(|\Lambda|^2), \\ 
L_{2} & = & u_{2} + O(|\Lambda|^2),  \\
\vdots&   &  \\
L_{8} & = & u_{8} + O(|\Lambda|^2), \\
L_{9} & = & u_{9}u_{10}+ O(|\Lambda|^3), \\
L_{10} & = & u_{9}u_{10} + O(|\Lambda|^3), \\
\end{array}
\right.
\end{equation*}
where $u_{9}=c_{0, 0, 2}$ and $u_{10}=-4c_{0, 1, 1} - 15c_{0, 2, 0} + 4c_{1, 0, 1} + 4b_{1, 1, 0}$. Therefore, we have 9 limit-cycles bifurcating from the center, being 8 with linear terms of Lyapunov constants and 1 with terms of order 2.

\subsection{Case $a_2=a_6=0$:}

For system \eqref{Valls} with $a_2=a_6=0$, the origin is a center if one of the following conditions holds:
\begin{itemize}
	\item[a)]$c_1=c_2=c_3=0$;
	\item[b)]$a_5=2c_1-c_2=c_3=0$; 
	\item[c)]$a_3=a_5=0$; 
\end{itemize}

In the item \textbf{a)}, assuming $a_1=1/2,a_3=-2,a_4=-2,a_5=7/8$, we have 4 limit-cycles using the linear terms of the Lyapunov constants.

For the item \textbf{b)}, taking $a_1=1/2,a_3=-2,a_4=-2,c_1=2/3$ we obtain the system 

\begin{equation*}\left\{\begin{array}{lll}
\dot{x} = y, \vspace{0.25cm} \\ 
\dot{y} = \dfrac{1}{2}x^2 - 2xz - 2y^2 - x, \vspace{0.25cm} \\
\dot{z} =\dfrac{2}{3}x^2 + \dfrac{4}{3}xy - z.
\end{array}
\right.
\end{equation*}
For this system, we obtain 9 limit-cycles, being 8 with the linear terms and 1 with the terms of order 2. The proof follows the same steps from the item \textbf{c)} in the case  $a_2=a_5=0$.

\subsection{Case $a_3=a_4=0$:}

For system \eqref{Valls} with $a_3=a_4=0$, the origin is a center if one of the following conditions holds:
\begin{itemize}
	\item[a)]$a_2=c_1=c_2=c_3=0$;
	\item[b)]$a_1=c_1=c_2=c_3=0$; 
	\item[c)]$a_6=2a_1a_2+a_5c_2=2c_1-c_2=c_3=0$; 
	\item[d)]$a_2=a_5=a_6=0$; 
	\item[e)]$a_1=a_5=a_6=0$; 
	\item[f)]$a_2=a_5=c_3=2c_1-c_2=0$; 
\end{itemize}

In the case  \textbf{a)}, taking $a_1=1/2,a_5=7/8,a_6=1$, we obtain 2 limit-cycles. Assuming $a_2=-2,a_5=7/8,a_6=1$ in the case  \textbf{b)} , the system obtained, when perturbed, present rank 3, that is, it is possible obtain 3 limit-cycles.

For the case \textbf{c)}, considering $a_2=-2,a_5=7/8,c_1=2/3$, we obtain the next system

\begin{equation*}\left\{\begin{array}{lll}
\dot{x} = y, \vspace{0.25cm} \\ 
\dot{y} =\dfrac{7}{24}x^2 - 2xy + \dfrac{7}{8}yz - x, \vspace{0.25cm} \\
\dot{z} =\dfrac{2}{3}x^2 + \dfrac{4}{3}xy - z.
\end{array}
\right.
\end{equation*}
Computing the Lyapunov constants up to order 3, we obtain 8 limit-cycles with the linear terms using the parameters perturbation $a_{0, 0, 2}$, $a_{0, 1, 1}, a_{0, 2, 0}, a_{1, 0, 1}, a_{1, 1, 0}, a_{2, 0, 0}, b_{0, 1, 1},c_{2, 0, 0}$, so we can write the Lyapunov constants as 

\begin{equation*}\left\{\begin{array}{lll}
L_{1} & = & u_{1} + O(|\Lambda|^2), \\ 
L_{2} & = & u_{2} + O(|\Lambda|^2),  \\
\vdots&   &  \\
L_{8} & = & u_{8} + O(|\Lambda|^2), \\
\end{array}
\right.
\end{equation*}
After an appropriate change of variable, we obtain the terms of order 2 of the Lyapunov constants $L_{9},L_{10}$ and $L_{11}$, i.e.:

\begin{equation*}\left\{\begin{array}{lll}
L_{9} & = & u_{9}u_{10} + O(|\Lambda|^3), \\ 
L_{10} & = & u_{9}u_{10} + O(|\Lambda|^3),  \\
L_{11} & = & u_{9}u_{10} + O(|\Lambda|^3), \\
\end{array}
\right.
\end{equation*}
where $u_{9}=992b_{0, 0, 2} - 147b_{0, 2, 0} + 336b_{1, 0, 1} - 504c_{0, 0, 2} - 147c_{1, 0, 1}$ and $u_{10}=992b_{0, 0, 2} - 147b_{0, 2, 0} + 336b_{1, 0, 1} - 434c_{0, 0, 2} - 147c_{1, 0, 1}$. So, we can obtain one more limit-cycle, adding up to 9 limit-cycles. Doing again an appropriate change of coordinates, we obtain the terms of order 3 of the Lyapunov constants. Then, we have

\begin{equation*}\left\{\begin{array}{lll}
L_{10} & = & u_{10}^{2}u_{11} + O(|\Lambda|^4),  \\
L_{11} & = & u_{10}^{2}u_{11} + O(|\Lambda|^4), \\
\end{array}
\right.
\end{equation*}
where $u_{11}=16905b_{0, 2, 0} - 3920b_{1, 0, 1} + 7595c_{0, 1, 1} + 22785c_{0, 2, 0} - 5880c_{1, 0, 1} - 288v_{10}$. Clearly we can obtain 1 more limit cycle. Therefore, we can obtain 10 limit-cycles with perturbations of this system.

In the case \textbf{d)}, considering $a1=1/2,c_3=2/8,c_1=2/3,c_2=-1/2$, we obtain 5 limit-cycles.
Analyzing case \textbf{e)}, assuming $a_2=-2,c_3=2/8,c_1:=2/3,c_2=-1/2$, we obtain 6 limit-cycles with the linear terms of the Lyapunov constants. For the case \textbf{f)}, taking $a_1=1/2,c_1=2/3$, we have 2 limit-cycles.

\subsection{Case $a_3=a_5=0$:}

For system \eqref{Valls} with $a_3=a_5=0$, the origin is a center if one of the following conditions holds:
\begin{itemize}
	\item[a)]$a_2=a_6=0$; 
	\item[b)]$a_2=2c_1-c_2=c_3=0$; 
	\item[c)]$a_1+a_4=a_6=0$; 
	\item[d)]$a_1+a_4=c_1=c_2=c_3=0$; 
\end{itemize}

Analyzing the item  \textbf{a)}, considering $a_1=1/2,a_4=-2,c_3=2/8,c_1=2/3,c_2=-1/2,$ we obtain the system 

\begin{equation*}\left\{\begin{array}{lll}
\dot{x} = y, \vspace{0.25cm} \\ 
\dot{y} =\dfrac{1}{2}x^2 - 2y^2 - x, \vspace{0.25cm} \\
\dot{z} =\dfrac{2}{3}x^2 - \dfrac{1}{2}xy + \dfrac{1}{4}y^2 - z.
\end{array}
\right.
\end{equation*}
Computing the Lyapunov constants up to order 3, we obtain 8 limit-cycles with the parameters perturbation $a_{0, 0, 2}$, $a_{0, 1, 1}$, $a_{0, 2, 0}$, $a_{1, 0, 1}$, $a_{2, 0, 0}$, $b_{0, 0, 2}$, $ b_{0, 1, 1}$, $b_{1, 0, 1}$. No more limit-cycles are obtained with order 2 and 3.

Studying the item  \textbf{b)}, taking $a_1=1/2,a_4=-2,a_6=1,c_1=2/3,$ we obtain the system 

\begin{equation*}\left\{\begin{array}{lll}
\dot{x} = y, \vspace{0.25cm} \\ 
\dot{y} =\dfrac{1}{2}x^2 - 2y^2 + z^2 - x, \vspace{0.25cm} \\
\dot{z} =\dfrac{2}{3}x^2 + \dfrac{4}{3}xy - z.
\end{array}
\right.
\end{equation*}
Computing the Lyapunov constants up to order 3, we obtain 8 limit-cycles with the parameters perturbation $a_{0, 0, 2}, a_{0, 1, 1}, a_{0, 2, 0}, a_{1, 0, 1}, a_{1, 1, 0}$, $a_{2, 0, 0}, b_{0, 1, 1},c_{2, 0, 0}$. No more limit-cycles is obtained with order 2 and 3.

For the item  \textbf{c)}, assuming $a_2=-2,a_4=-2,c_3=2/8,c_1=2/3,c_2=-1/2$, we obtain 7 limit-cycles. In the item \textbf{d)}, taking $a_2=-2,a_4=-2$ we obtain 2 limit-cycles.

\subsection{Case $a_3=a_6=0$:}

For system \eqref{Valls} with $a_3=a_6=0$, the origin is a center if one of the following conditions holds:
\begin{itemize}
	\item[a)]$a_2=a_5=0$; 
	\item[b)]$a_5=a_1+a_4=0$; 
	\item[c)]$a_2=c_1=c_2=c_3=0$; 
	\item[d)]$a_1+a_4=c_1=c_2=c_3=0$; 
	\item[e)]$2a_1a_2+a_5c_2=a_4=2c_1-c_2=c_3=0$; 
\end{itemize}

In the item \textbf{a)}, taking $a_1=1/2,a_4=-2,c_3=2/8,c_1=2/3,c_2=-1/2$, we obtain the system

\begin{equation*}\left\{\begin{array}{lll}
\dot{x} = y, \vspace{0.25cm} \\ 
\dot{y} =\dfrac{1}{2}x^2 - 2y^2- x, \vspace{0.25cm} \\
\dot{z} =\dfrac{2}{3}x^2 - \dfrac{1}{2}xy + \dfrac{1}{4}y^2 - z.
\end{array}
\right.
\end{equation*}
Computing the Lyapunov constants up to order 3, we obtain 8 limit-cycles with the linear terms of Lyapunov constants $a_{0, 0, 2}, a_{0, 1, 1}, a_{0, 2, 0}$, $a_{1, 0, 1}, a_{2, 0, 0}, b_{0, 0, 2}, b_{0, 1, 1},b_{1, 0, 1}.$ No more limit-cycles is obtained with Lyapunov constants up to order 2 and 3.

Analyzing item \textbf{b)} and considering $a_2=-2,a_4=-2,c_3=-2/8,c_1=4/3,c_2=-1/2$, we obtain 7 limit-cycles with linear terms of Lyapunov constants.

For the item \textbf{c)} taking $a_1=1/2, a_4=-2,a_5=7/8$ we obtain 4 limit-cycles. In the item \textbf{d)}, assuming $a_2=-2,a_4=-2,a_5=7/8$, again we obtain 4 limit-cycles.

In the item \textbf{e)}, assuming $a_2=-1,a_5=-8/9,c_2=-1/2$, we have

\begin{equation*}\left\{\begin{array}{lll}
\dot{x} = y, \vspace{0.25cm} \\ 
\dot{y} =\dfrac{2}{9}x^2 -xy- \dfrac{8}{9}yz-x, \vspace{0.25cm} \\
\dot{z} =-\dfrac{1}{4}x^2 - \dfrac{1}{2}xz - z.
\end{array}
\right.
\end{equation*}
The linear terms has rank 8 with the perturbation terms $a_{0, 0, 2}, a_{0, 1, 1}, a_{0, 2, 0}$, $a_{1, 0, 1}, a_{1, 1, 0}, a_{2, 0, 0}, b_{0, 0, 2}, c_{0, 2, 0}$. Making an appropriate change of coordinates, we obtain

\begin{equation*}\left\{\begin{array}{lll}
L_9  = u_{9}u_{10}+O(|\Lambda|^3) \vspace{0.25cm} \\ 
L_{10} = u_{9}u_{10}+O(|\Lambda|^3),
\end{array}
\right.
\end{equation*}
where $u_{9}=3159b_{0, 1, 1} + 11772b_{0, 2, 0} - 972b_{1, 0, 1} - 2808b_{1, 1, 0} + 12636b_{2, 0, 0} - 864c_{1, 0, 1} + 2912c_{1, 1, 0} + 5408c_{2, 0, 0}$ and $u_{10}=3159b_{0, 1, 1} + 11772b_{0, 2, 0} - 972b_{1, 0, 1} - 2808b_{1, 1, 0} + 12636b_{2, 0, 0} + 1458c_{0, 0, 2} - 864c_{1, 0, 1} + 2912c_{1, 1, 0} + 5408c_{2, 0, 0}$. Therefore, 9 limit-cycles unfold from the center.

\subsection{Case $a_4=a_5=0$:}

For system \eqref{Valls} with $a_4=a_5=0$, the origin is a center if one of the following conditions holds:
\begin{itemize}
	\item[a)]$a_1=c_1=c_2=c_3=0$; 
	\item[b)]$a_1=a_6=2c_1-c_2=c_3=0$; 
	\item[c)]$a_2=2c_1-c_2=c_3=0$; 
	\item[d)]$a_2=a_3=a_6=0$; 
	\item[e)]$a_1=a_3=a_6=0$; 
\end{itemize}

Studying the item \textbf{a)}, assuming $a_2=-1,a_3=-2,a_6=1$, we obtain 3 limit-cycles. In the item \textbf{b)}, considering $a_2=-1,a_3=-2,c_1=4/3$, we obtain 5 limit-cycles. For the item \textbf{c)}, taking $a_1=1/2,a_3=-2,a_6=1,c_1=4/3$, we obtain 7 limit-cycles. In the item \textbf{d)}, considering $a_1=1/2,c_3=-2/8,c_1=4/3,c_2=-1/2$ we obtain 5 limit-cycles. Finally, in item \textbf{e)}, taking $a_2=-1,c_3=-2/8,c_1=4/3,c_2=-1/2$, we have 6 limit-cycles.

\subsection{Case $a_4=a_6=0$:}

For system \eqref{Valls} with $a_4=a_6=0$, the origin is a center if one of the following conditions holds:
\begin{itemize}
	\item[a)] $a_2=c_1=c_2=c_3=0$;
	\item[b)] $a_1=c_1=c_2=c_3=0$;
	\item[c)] $a_2=a_5=2c_1-c_2=c_3=0$;
	\item[d)] $a_1=a_5=2c_1-c_2=c_3=0$;
	\item[e)] $a_3=2a_1a_2-a_5c_2=2c_1-c_2=c_3=0$;
	\item[f)] $a_2=a_3=a_5=0$;
	\item[g)] $a_1=a_3=a_5=0$;
\end{itemize}

In the cases \textbf{a)} and \textbf{b)}, taking $a_1=1/2,a_3=-2,a_5=-8/9$ and $
a_2=-1,a_3=-2,a_5=-8/9$ respectively, we obtain 3 limit-cycles in each case. For the cases \textbf{c)} and \textbf{d)}, considering $a_1=1/2,a_3=-2,c_1=2/3$ and $a_2=-1,a_3=-2,c_1=-2/3$ respectively, we have 5 limit-cycles in each case.

Studying item \textbf{e)} assuming $a_5=1$, $a_2=-3/2$, $c_2=-1/2$, we obtain the system

\begin{equation*}\left\{\begin{array}{lll}
\dot{x} = y, \vspace{0.25cm} \\ 
\dot{y} =-\dfrac{1}{6}x^2 - \dfrac{3}{2}xy + yz - x, \vspace{0.25cm} \\
\dot{z} =-\dfrac{1}{4}x^2 - \dfrac{1}{2}xy - z.
\end{array}
\right.
\end{equation*}
Their rank is 8 with the perturbation parameters $a_{0, 0, 2}, a_{0, 1, 1}, a_{0, 2, 0}, a_{1, 0, 1}$, $a_{1, 1, 0}, a_{2, 0, 0}, b_{0, 0, 2}, b_{0, 1, 1},c_{2, 0, 0}$. Making an appropriate change of coordinates to vanish the linear terms of $L_{9}$ and $L_{10}$, we obtain 

\begin{equation*}\left\{\begin{array}{lll}
L_{9} & = & u_{9}u_{10} + O(|\Lambda|^3),  \\
L_{10} & = & u_{9}u_{10} + O(|\Lambda|^3), \\
\end{array}
\right.
\end{equation*}
where $u_{9}=42a_{2, 0, 0} + 6b_{0, 2, 0} - 9b_{1, 0, 1} + 14c_{0, 2, 0} + 6c_{1, 0, 1}$ and $u_{10}=168a_{2, 0, 0} + 24b_{0, 2, 0} - 36b_{1, 0, 1} + 243c_{0, 0, 2} + 56c_{0, 2, 0} + 24c_{1, 0, 1}$. No more limit-cycle is obtained with order three. Therefore, at least 9 limit-cycles unfold from the center.

For the cases \textbf{f)} and \textbf{g)}, assuming $a_1=1/2,c_3=-2/5,c_1=-2/3,c_2=-1/2$ and $a_2=-1,c_3=-2/5,c_1=-2/3,c_2=-1/2$ respectively, we obtain 5 and 6 limit-cycles respectively.

\subsection{Case $a_5=a_6=0$:}

For system \eqref{Valls} with $a_5=a_6=0$, the origin is a center if one of the following conditions holds:
\begin{itemize}
	\item[a)] $a_2=a_3=0$;
	\item[b)] $a_1+a_4=a_3=0$;
	\item[c)] $a_1+a_4=c_1=c_2=c_3=0$;
	\item[d)] $a_2=2c_1-c_2=c_3=0$;
	\item[e)] $a_1=a_4=2c_1-c_2=c_3=0$.
\end{itemize}

For the item \textbf{a)}, considering $a_1=1,a_4=-2,c_3=-2/5,c_1=-2/3,c_2=2$, we have the system 

\begin{equation*}\left\{\begin{array}{lll}
\dot{x} = y, \vspace{0.25cm} \\ 
\dot{y} =x^2 - 2y^2 - x, \vspace{0.25cm} \\
\dot{z} =-\dfrac{2}{3}x^2 + 2xy - \dfrac{2}{5}y^2 - z.
\end{array}
\right.
\end{equation*}
For this system, we have 8 limit-cycles using the linear terms of Lyapunov constants with the perturbation parameters $a_{0, 0, 2}, a_{0, 1, 1}, a_{0, 2, 0}$, $a_{1, 0, 1}$, $a_{2, 0, 0}, b_{0, 0, 2},b_{0, 1, 1},b_{1, 0, 1}.$ No more limit-cycles are obtained with Lyapunov constants up to order 2 and 3.

For the case \textbf{b)}, considering $a_2=-1,a_4=-2,c_3=65,c_1=2/56,c_2=1$, we obtain the system

\begin{equation*}\left\{\begin{array}{lll}
\dot{x} = y, \vspace{0.25cm} \\ 
\dot{y} =2x^2 - xy - 2y^2 - x, \vspace{0.25cm} \\
\dot{z} =\dfrac{1}{28}x^2 + xy + 65y^2 - z.
\end{array}
\right.
\end{equation*}
Perturbing this system, we obtain 8 limit-cycles using the linear terms of Lyapunov constants with the perturbation parameters $a_{0, 0, 2}$, $a_{0, 1, 1}$, $a_{0, 2, 0}$, $a_{1, 0, 1}$, $b_{0, 0, 2}$, $b_{0, 1, 1}$, $b_{0, 2, 0}$, $b_{1, 0, 1}$. No more limit-cycles are obtained with the Lyapunov constants up to order 2 and 3.

Analyzing item \textbf{c)}, assuming $a_2=-1,a_4=-2,a_3=-8/9$, we obtain 4 limit-cycles. 

For the case \textbf{d)}, considering $a_1=-1,a_4=2,a_3=-9,c_1=2/5$, we obtain the system 

\begin{equation*}\left\{\begin{array}{lll}
\dot{x} = y, \vspace{0.25cm} \\ 
\dot{y} =-x^2 - 9xz + 2y^2 - x, \vspace{0.25cm} \\
\dot{z} =\dfrac{2}{5}x^2 + \dfrac{4}{5}xy - z.
\end{array}
\right.
\end{equation*}
Computing the Lyapunov constants up to order 3, we have 8 limit-cycles using the linear terms of the Lyapunov constants with the perturbation parameters $a_{0, 0, 2}, a_{0, 1, 1}, a_{0, 2, 0},$ $a_{1, 0, 1}, a_{1, 1, 0}, a_{2, 0, 0}, b_{0, 1, 1},c_{2, 0, 0}$. So, we can write the first 8 Lyapunov constants as follows

\begin{equation*}\left\{\begin{array}{lll}
L_{1} & = & u_{1} + O(|\Lambda|^2), \\ 
L_{2} & = & u_{2} + O(|\Lambda|^2),  \\
\vdots&   &  \\
L_{8} & = & u_{8} + O(|\Lambda|^2). \\
\end{array}
\right.
\end{equation*}
Doing an appropriate change of coordinates for vanishing the linear terms of the $L_{9}$ and $L_{10}$, we obtain
\begin{equation*}\left\{\begin{array}{lll}
L_{9} & = & u_{9}u_{10} + O(|\Lambda|^3), \\ 
L_{10} & = & u_{9}u_{10} + O(|\Lambda|^3),  \\
\end{array}
\right.
\end{equation*}
where $u_{9}=c_{0, 0, 2}$ and $u_{10}=2c_{1, 0, 1} + 15c_{0, 2, 0} - 2c_{0, 1, 1} + 2b_{1, 1, 0}$. Then, it is clearly that we can obtain one more limit-cycle. No more limit-cycles are obtained with Lyapunov constants up to order 3. Therefore, we can obtain 9 limit-cycles.

In the item \textbf{e)}, taking $a_2=3,a_3=5,c_1=2$ we have 5 limit-cycles.

\section{Edneral-Mahdi-Romanovski-Shafer quadratic systems}

In the paper \cite{Adam} the authors extend the Lyapunov method to solve the center problem for the following family of quadratic three-dimensional systems having a Hopf singularity.
\begin{equation*}\left\{\begin{array}{lll}
\dot{x}&=&-y+ax^2+ay^2+cxz+dyz, \vspace{0.25cm} \\ 
\dot{y}&=&x+bx^2+by^2+exz+fyz, \vspace{0.25cm} \\
\dot{z}&=&-z+Sx^2+Sy^2+Txz+Uyz.
\end{array}
\right.
\end{equation*}
We go through every branch of the Bautin Variety of system \eqref{Adam} and compute the Lyapunov constants of the perturbation \eqref{eqPerturbation}.

\subsection{Branch $a=b=c+f=0$, $S=1$, with $c=d+e=0$:} We compute the first 11 Lyapunov constants for a generic point $e=1,T=\frac{1}{2},U=\frac{3}{4}$ on this branch of the Bautin Variety. The system is given by 

\begin{equation*}\left\{\begin{array}{lll}
\dot{x}&=& -yz - y, \vspace{0.25cm} \\ 
\dot{y}&=& xz + x, \vspace{0.25cm} \\
\dot{z}&=& x^2 + \frac{1}{2}xz + y^2 + \frac{3}{4}yz - z.
\end{array}
\right.
\end{equation*}
Their rank is $9$, and by Theorem \ref{TeoBifLinear} the cyclicity is at least $9$. However considering terms of order 2, 3 and 4 of the Lyapunov constants yields no additional limit-cycles.

\subsection{Branch $a=b=c+f=0$, $S=1$, with $8c+T^2-U^2=4(e-d)-T^2-U^2=2(e+d)+TU=0$:}
We compute the first 11 Lyapunov constants for a generic point $e=0,U=1,T=1$ on this branch of the Bautin Variety. The respective system is given by

\begin{equation*}\left\{\begin{array}{lll}
\dot{x}&=& -\frac{1}{2}yz - y, \vspace{0.25cm} \\ 
\dot{y}&=& x, \vspace{0.25cm} \\
\dot{z}&=& x^2 + xz + y^2 + yz - z.
\end{array}
\right.
\end{equation*}
Their rank is $9$, and by Theorem \ref{TeoBifLinear} the cyclicity is at least $9$. Considering the Lyapunov constants up to order 2, 3 and 4 yields no additional limit-cycles.

\subsection{Branch $d+e=c=f=0$, $S=1$, with $a=b=0$:}
We compute the first 11 Lyapunov constants for a generic point $e=-2,T=\frac{1}{2},U=\frac{3}{4}$ on this branch of the Bautin Variety. Then, the system is given by

\begin{equation*}\left\{\begin{array}{lll}
\dot{x}&=& 2yz - y, \vspace{0.25cm} \\ 
\dot{y}&=& -2xz+x, \vspace{0.25cm} \\
\dot{z}&=& x^2 + \frac{1}{2}xz + y^2 + \frac{3}{4}yz - z.
\end{array}
\right.
\end{equation*}
Their rank is $9$, and by Theorem \ref{TeoBifLinear} the cyclicity is at least $9$. No more limit-cycles are obtained with Lyapunov constants up to order 2, 3 and 4.

\subsection{Branch $d+e=c=f=0$, $S=1$, with $T-2a=U-2b=0$:}
We compute the first 11 Lyapunov constants for a generic point $a=1,b=\frac{5}{2},e=-2$ on this branch of the Bautin Variety. 

\begin{equation*}\left\{\begin{array}{lll}
\dot{x}&=& x^2 + y^2 + 2yz - y, \vspace{0.25cm} \\ 
\dot{y}&=& \frac{5}{2}x^2 - 2xz + \frac{5}{2}y^2 + x, \vspace{0.25cm} \\
\dot{z}&=& x^2 + 2xz + y^2 + 5yz - z.
\end{array}
\right.
\end{equation*}
Their rank is $9$, and by Theorem \ref{TeoBifLinear} the cyclicity is at least $9$. No additional limit-cycles are obtained with Lyapunov constants up to order 2, 3 and 4.

\subsection{Branch $d+e=c=f=0$, $S=1$, with $d=e=0$:}
We compute the first 11 Lyapunov constants for a generic point $a=1,b=\frac{5}{2},T=\frac{3}{8},U=-4$ on this branch of the Bautin Variety. 

\begin{equation*}\left\{\begin{array}{lll}
\dot{x}&=& x^2 + y^2 - y, \vspace{0.25cm} \\ 
\dot{y}&=& \frac{5}{2}x^2 + \frac{5}{2}y^2 + x, \vspace{0.25cm} \\
\dot{z}&=& x^2 + \frac{3}{8}xz + y^2 - 4yz - z.
\end{array}
\right.
\end{equation*}
Their rank is $8$, and by Theorem \ref{TeoBifLinear} the cyclicity is at least $8$. Considering the Lyapunov constants up to order 2, 3 and 4 yields no additional limit-cycles.

\subsection{Branch $S=0$:}
We compute the first 11 Lyapunov constants for a generic point $a=1,b=2,c=3,d=-\frac{1}{2},e=\frac{5}{2},f=\frac{2}{3},T=-1,U=2$ on this branch of the Bautin Variety. 

\begin{equation*}\left\{\begin{array}{lll}
\dot{x}&=& x^2 + 3xz + y^2 - \frac{1}{2}yz - y, \vspace{0.25cm} \\ 
\dot{y}&=& 2x^2 + \frac{5}{2}xz + 2y^2 + \frac{2}{3}yz + x, \vspace{0.25cm} \\
\dot{z}&=& -xz + 2yz - z.
\end{array}
\right.
\end{equation*}
Their rank is $5$, and by Theorem \ref{TeoBifLinear} the cyclicity is at least $5$. Studying the Lyapunov constants up to order 2 and 3 yields two additional limit-cycles. Thus, the cyclicity is at least $7$.

\section{Final Comments}

In our work, we made a chart of the cyclicity of quadratic systems having a Hopf singular points and the computations suggests that it is very difficult to obtain examples of centers bifurcating more than 10 limit-cycles. In the literature, besides the systems we worked in the present paper, \cite{PeiYu} present one more example for which 10 limit-cycles can unfold. For 11 limit-cycles, the only known example is presented in \cite{Ivan}. The authors of \cite{Ivan} also conjectured that the maximum number of limit-cycles bifurcating from centers of systems \eqref{eq1} is 12 and, from the extensive list of systems that were considered in our work, we found only one system that reaches this bound.

\section{Acknowledgments}

We would like to thank professors Joan Torregrosa and Claudio Pessoa for their immensely appreciated support, helpful discussions and valuable suggestions which surely made the present paper enriched. The first author is supported by S\~ao Paulo Research Foundation (FAPESP) grant 20/04717-0. The second author is supported by S\~ao Paulo Research Foundation (FAPESP) grant 19/13040-7.

\addcontentsline{toc}{chapter}{Bibliografia}
\bibliographystyle{siam}
\bibliography{Referencias_HopfR3.bib}
\end{document}